%%%%%%%%%%%%%%%%%%%%%%%%%%%%%%AMS-TeX fil%%%%%%%%%%%%%%%%%%%%%%%%%
\magnification=\magstep1
%%%%%%%%%%%%%%%%%%%%%%%A4?%%%%%%%%%%%%%%%%%%%%%%%%%%%%%
\hsize=6.1 true in
\vsize=9.3 true in
\hoffset=0pt
\voffset=0pt
%\hcorrection{2,5truemm}
%\topskip=2cm
%\voffset -0.4in
%\vsize 24truecm
%\hsize 15.5truecm

\input amstex
\documentstyle{amsppt}
%%%%%%%%%%%%%%%%%%%%%%%%spesialiteter%%%%%%%%%%%%%%%%%%%%%%%%%%%%%
\TagsOnRight
\NoRunningHeads
%\nologo
%%%%%%%%%%%%%%%%%%%%%%%makro for datum%%%%%%%%%%%%%%%%%%%%%%%%%%%%%%%
  \def\today{\ifcase\month \or January \or February \or March \or April
    \or May \or June \or July \or August \or September \or October
    \or November \or December \fi \number \day, \number \year}
%%%%%%%%%%%%%%%%%%%%%makro for notat i margen%%%%%%%%%%%%%%%%%%%%%%
\newdimen\notespace \notespace=5mm %avstand fra venstre marg
\newdimen\maxnote  \maxnote=13mm

\def\obs{\strut\vadjust{\kern-\dp\strutbox
  \vtop to\dp\strutbox{\vss \baselineskip=\dp\strutbox
  \moveleft\notespace\llap{\hbox to \maxnote{\hfil$!$}}\null}}\ }

\def\cref#1{{\rm #1}}
\def\pref#1{{\rm(#1)}}
%\def\pref#1{{\rm(\marg#1)}}

%\redefine\marg{\relax\ }

%%%%%%%%%%%%%%%%%%diverse makroer for format%%%%%%%%%%%%%%%%%%%%%
\def\risom{\,\overset{\sim}\to{\smash\longrightarrow\vrule
    height0.5ex width0pt depth0pt}\,}
\def\dfn#1{{\it #1}}
%%%%%%%%%%%%%%%%%%diverse makroer for letthet%%%%%%%%%%%%%%%%%%
\font\mibold=cmmib10
\font\miibold=cmmib8
\define\boldkappa{\text{\mibold \char'024}}
\define\kap#1{\boldkappa(#1)}

\define\doublearrow#1#2{\overset #1\to{\underset #2\to
\rightrightarrows}}

\define\F#1{{\Bbb F}_{\!\!#1}}

\def\-alg{\operatorname{\kern-1.0pt{-alg}}}

\def\alg{\operatorname{alg}}

\def\dual{\,\check{}\,}
\def\End{\operatorname{End}}
\def\ev{\operatorname{ev}}

\def\id{\operatorname{id}}
\def\Hilb{\Cal Hilb}
\def\Hom{\operatorname{Hom}}

\def\polring{A[X_1, X_2, \dots, X_m]}

\def\Proj{\operatorname{Proj}}

\def\Spec{\operatorname{Spec}}

\def\Sym{\operatorname{Sym}}

\redefine\det{\operatorname{det}}

%%%%%%%%%%%%%%%%%%%%%%%%%%%%hyphenation%%%%%%%%%%%%%%%%%%%%%%%%
\hyphenation{homo-mor-phism}

%%%%%%%%%%%%%%%%%%%%%%%topmatter%%%%%%%%%%%%%%%%%%%%%%%%%%%%%%%%%%%
\topmatter
\title An elementary, explicit, proof of the existence of
Hilbert schemes of points
\endtitle   
\author T. S. Gustavsen, D. Laksov,  R. M. Skjelnes
\endauthor
\affil Department of Mathematics, University of Oslo,  Department of
Mathematics, KTH 
\endaffil 
\address   NO-0316 Oslo-Blindern, Norway, S-100 44 Stockholm, Sweden 
\endaddress 
\email stolen\@math.uio.no, laksov\@math.kth.se, skjelnes\@math.kth.se
\endemail 
%\date \today\ ({\tt \jobname})\enddate
\abstract

 We give a short, elementary, and explicit proof of the existence of
 the Hilbert scheme of $n$ points of an affine scheme $\Spec(R)$ over
 an affine base scheme $\Spec(A)$. As a consequence we obtain an easy
 proof of the existence of the Hilbert scheme of $n$ points of a
 projective scheme $X$ over an arbitrary base scheme $S$, that is,
 when $X=\Proj(\Cal R)\to S$ where $\Cal R$ is any quasi-coherent
 graded $\Cal O_S$-algebra.

 Our methods rely on simple algebraic constructions and avoid the
 usual embeddings into high dimensional spaces, and thus, the use of
 Castelnuovo-Mumford regularity. We give explicit expressions of
 members of an affine covering of the Hilbert schemes, involving few
 variables satisfying natural equations.
\endabstract

\keywords Hilbert scheme of points, endomorphism ring, evaluation map, commuting matrices  \endkeywords

\subjclass Primary 14C05. Secondary 15A27, 16S50 \endsubjclass

%\toc
%\widestnumber\subhead{1.2}
%\head {\phantom{0}} Introduction\endhead
%\head 1. Coverings of the Hilbert functor\endhead
%\head 2. The Hilbert functor and endomorphisms\endhead 
%\head 3. Representing endomorphisms\endhead
%\head 4. Representing commuting endomorphisms\endhead
%\head 5. Sections and closed subschemes\endhead
%\endtoc

\endtopmatter

%%%%%%%%%%%%%%%%%%%%%dokument%%%%%%%%%%%%%%%%%%%%%%%%%%
 \document

 \head Introduction\endhead

 The main objective of this note is to give a short, elementary, and
 explicit proof of the existence of the Hilbert scheme of $n$ points
 of an affine scheme $\Spec(R)$ over an affine base scheme $\Spec(A)$,
 where $R$ is any $A$-algebra. As a consequence of our results we
 obtain an easy proof of the existence of the Hilbert scheme of $n$
 points of a projective scheme $X$ over an arbitrary base scheme $S$,
 that is, when $X=\Proj(\Cal R)\to S$ where $\Cal R$ is any
 quasi-coherent  graded $\Cal O_S$-algebra.

When $S$ is locally noetherian and $\Cal R$ is locally finitely
generated by elements of degree one, our result follows from the more
general results of Grothendieck
\cite{G} (see also \cite{A}, \cite{H}, \cite{HS}, \cite{M}, \cite{N},
or \cite{St}). However, our methods rely only on simple algebraic
constructions and avoid embeddings into high dimensional spaces via
Castelnuovo-Mumford regularity. We give explicit expressions for
members of an affine covering of the Hilbert schemes, involving few
variables satisfying natural equations. This provides a powerful tool
for studying the geometric properties of Hilbert schemes of points. In
particular we obtain a natural description of an open subset of the
\dfn{generic} component of the Hilbert schemes, corresponding to $n$
distinct points.

One of the main features of our construction is that it gives explicit
equations defining the Hilbert schemes of $n$ points of $\Spec(R)$ as
a closed subscheme of the scheme of commuting $n\times n$-matrices.

Furthermore, to illustrate the flexibility of our methods we compute
the Hilbert scheme of the scheme $\Spec (S^{-1}A[X])$ over $\Spec(A)$,
where $S$ is a multiplicatively closed subscheme of the polynomial
ring $A[X]$ in the variable $X$ over $A$ (see \cite{LS}, \cite{LST},
and \cite{S}). 

Our construction is based upon two simple ideas. The first is to use
the well-known description of the $R$-module structure on an
$A$-module $F$ in terms of $A$-algebra homomorphisms from $R$ to
$\End_A(F)$. The second is the observation that the $A$-module
homomorphisms from an $A$-module $M$ to $\End_A(F)$ correspond to
$A$-algebra homomorphisms from the symmetric algebra
$\Sym_A(M\otimes_A\End_A(F)\dual)$ to $A$. From the latter observation
we obtain that the $A$-algebra homomorphisms from $\Sym_A(M)$ to
$\End_A(F)$ correspond to the algebra homomorphisms $H\to A$, where
$H$ is the residue of the algebra $\Sym_A(M\otimes_A\End_A(F)\dual)$
by the ideal corresponding to commuting $n\times n$-matrices, and with
coefficients in $\Sym_A(M\otimes_A\End_A(F)\dual)$.

 \def\sectno{1}

 \head \sectno. The basic functors\endhead

 In this section we define the basic functors used in the remaining
 part of the article, and we give an explicit relation between
 $A$-module homomorphisms $R\to F$ from an $A$-algebra $R$ to an
 $A$-module $F$, and the $A$-algebra homomorphisms $R\to \End_A(F)$ to
 the endomorphisms of $F$.

 \definition{\sectno.1 Notation}Let $A$ be a commutative ring with
 unit, and let $F$ be a free $A$-module of finite rank with a
 distinguished element $e$ that is part of a basis of $F$.
 Furthermore we let $R$ be an $A$-algebra. We write
 $$\End_A(F) @>{\ev_e}>>F$$
  for the \dfn{evaluation map} defined by $\ev_e(u) =u(e)$.

 We shall be using several functors from algebras to sets that we
 shall define next. 
 \enddefinition

 \definition{\sectno.2 Functors of homomorphisms of modules} 
 For every pair of $A$-modules $M$ and $N$ we denote by $\Cal
 Hom(M,N)$ the functor from $A$-algebras to $A$-modules that on an
 $A$-algebra $B$ takes the value $\Cal Hom_B(M,N)
 =\Hom_B(B\otimes_AM, B\otimes_AN)$. When $M=N$ we write $\Cal
 End(M)=\Cal Hom(M,M)$. For every $A$-algebra homomorphism
 $\varphi:B\to C$ we denote the natural homomorphism by 
 $$\Cal Hom_\varphi: \Cal Hom_B(M,N)\to \Cal Hom_C(M,N).$$
 \enddefinition

 \definition{\sectno.3 Functors of homomorphisms of algebras}
  Denote by $\Cal Hom_{\alg}(R, \End(N))$ the functor from
  $A$-algebras to sets that on an $A$-algebra $B$ takes the value
 $$\Cal Hom_{B\-alg}(R,\End(N))
 =\Hom_{B\-alg}(B\otimes_AR,B\otimes_A\End(N)).$$

 Moreover, denote by $\Cal H_R$ the functor from $A$-algebras to
 sets that on the $A$-algebra $B$ takes the value $\Cal H_R(B)
 =\Hom_{A\-alg}(R,B)$. From the canonical isomorphism
 $\Hom_{A\-alg}(R,A) \to \Hom_{B\-alg}(B\otimes_AR, B\otimes_AA)$ we
 obtain a natural isomorphism of functors
 $$\Cal H_R \to \Cal Hom_{\alg}(R,A).$$
 \enddefinition

 \definition{\sectno.4 Subfunctors defined by a unit} Recall that we
 have a distinguished element $e$ of our free $A$-module $F$. For any
 $A$-algebra $B$ we let ${\Cal Hom}^e(R,F)$ consist of surjective
 $B$-module homomorphisms $u : B\otimes_A R \to B\otimes_A F$ such
 that the kernel is an ideal in $B\otimes_AR$ and where $u(1_B\otimes
 1_R)=1_B\otimes e$. We have that ${\Cal Hom}^e(R,F)$ is a subfunctor
 of ${\Cal Hom}(R,F)$.  \enddefinition
 
 \definition{\sectno.5 Subfunctors defined by sections}
 For every $A$-module homomorphism 
 $$\beta :F\to R \quad\text{such that}\quad \beta(e) =1_R$$
 we denote by $\Cal Hom^\beta(R,F)$ the subfunctor of $\Cal
 Hom^e(R,F)$ whose value on an $A$-algebra $B$ that consists of
 $B$-module homomorphisms $u:B\otimes_AR \to 
 B\otimes_AF$ such that the composite homomorphism 
  $$B\otimes_AF @>{\id_B\otimes\beta}>> B\otimes_AR @>{u}>>
 B\otimes_AF$$
 is the identity.

 Denote by $\Cal Hom_{\alg}^\beta(R,\End(F))$ the subfunctor of
 $\Cal Hom_{\alg}(R,\End(F))$ that on the $A$-algebra $B$ consists of
 the homomorphisms lying in $\Cal Hom_B^\beta(R,\End(F))$.
 \enddefinition

 \definition{\sectno.6 Remark}We note that for an $A$-module $N$ there
 is a natural structure as a left $\End_A(N)$-module on $N$, and that
 there is a natural correspondence between $R$-module structures on
 $N$ and homomorphisms of $A$-algebras $R\to \End_A(N)$.
\enddefinition

 \proclaim{\sectno.7 Proposition} Let $\Cal H\Cal E_R^e$ be the
  functor from $A$-algebras to sets whose value $\Cal H\Cal E_R^e(B)$
  at the $A$-algebra $B$ consists of the $B$-algebra homomorphisms
  $\varphi:R\to \End_B(B\otimes_AF)$ such that the composition of the
  homomorphisms  $B\otimes_AR
  @>{\varphi}>> \End_B(B\otimes_AF) @>{\ev_{1_B\otimes e}}>>
  B\otimes_AF$ is surjective.

  There is a natural isomorphism of functors
 $$\Cal H\Cal E_R^e \to \Cal Hom^e(R,F)\tag{\sectno.7.1}$$
 that, for every $A$-algebra $B$, maps a $B$-algebra homomorphism
 $\varphi:B\otimes_AR\to \End(B\otimes_AF)$ to $\ev_{1_B\otimes
 e}\varphi:B\otimes_AR\to B\otimes_AF$. The inverse of the latter
 homomorphism  maps a $B$-module
 homomorphism $u:B\otimes_AR\to B\otimes_AF$ to the $B$-algebra
 homomorphism $\varphi:B\otimes_AR\to \End_B(B\otimes_AF)$ that makes
 $B\otimes_AF$ into the unique $B$-algebra such that $u$ becomes a
 $B$-algebra homomorphism.

  The isomorphism \pref{\sectno.7.1}  induces a natural isomorphism of
 functors  
$$\Cal H\Cal E_R^\beta \to \Cal Hom^\beta(R,F),$$ where $\Cal H\Cal
 E_R^\beta$ is the subfunctor of $\Cal H\Cal E_R^e$ consisting of the
 elements $\varphi\in\Cal H\Cal E_R(B)$ such that the composite $
 B\otimes_AF @>{\id_B\otimes_A\beta}>> B\otimes_A R @>{\varphi}>>
 \End_B(B\otimes_A F) @>{\ev_{1_B\otimes e}}>> B\otimes_A F$ is the
 identity.  \endproclaim

  \demo{Proof}The map \pref{\sectno.7.1} described in the proposition
is clearly functorial.

 To prove the first part of the proposition it thererfore clearly
 suffices to show the proposition in the case $A=B$. Let $\varphi:R\to
 \End_A(F)$ be a homomorphism in $\Cal H\Cal E_R^e(A)$ and let
 $u=\ev_e\varphi:R\to F$. Then $u$ is surjective by 
 definition and $u(1_R) =\ev_e\varphi(1_R) =\ev_e(\id_{\End_A(F)})
 =\id_{\End_A(F)}(e) =e$. Moreover, the kernel of $\ev_e$ is a left
 ideal in $\End_A(F)$ and consequently the kernel of $u$ is an ideal
 in $R$. We have thus constructed a map from $\Cal H\Cal E_R^e(A)$
 to $\Cal  Hom_A^e(R,F)$.

 Conversely, let $u:R\to F$ be in $\Cal Hom_A^e(R,F)$. Then $F$ has a
 unique $A$-algebra structure such that $u$ is a homomorphism of
 algebras. Let
 $$\varphi:R\to \End_A(F)$$
 be the homomorphism such that the image $\varphi_f$ of $f$ is defined
 by $\varphi_f(x) =u(f)x$, where the product on the right side is
 multiplication in $F$ with the given algebra structure. It is clear
 that $\varphi$ is an $A$-algebra homomorphism. We have
 $\ev_e\varphi_f 
 =\varphi_f(e) =u(f)e =u(f)$. Hence $u=\ev_e\varphi$. In particular
 $\ev_e\varphi$ is surjective so that $\varphi$ is in $\Cal
 H_A^e(R,\End(F))$. We have thus constructed a map from $\Cal
 Hom_A^e(R,F)$ to $\Cal H\Cal E_R^e(A)$. 

 It is clear that the two maps that we have constructed are inverses.
 
 The last part of the proposition is also clear.
 \enddemo

 The following result is often useful to compute representants of
 functors. 

\proclaim{\sectno.8 Proposition}Let $\varphi:R\to \End_A(F)$ be an
 $A$-algebra homomorphism such that $\ev_e\varphi$ is surjective, and
 give $F$ the corresponding structure as an $A$-algebra such that
 $\ev_e\varphi$ is an $A$-algebra homomorphism. Then an endomorphism
 $v\in \End_A(F)$ commutes with all elements in $\varphi(R)$ if and
 only if 
 $$v(x)=v(e)x$$
 for all $x\in F$.  In other words, if
$$\chi_\varphi:F\to \End_A(F)$$
 is the $A$-algebra homomorphism defined by $\chi_\varphi(y)(x)=yx$
 for all $x$ and $y$ in $F$, then $\chi_\varphi(F)$ is the subset of
 $\End_A(F)$ of elements that commute with all elements in
 $\varphi(R)$.
 \endproclaim

 \demo{Proof}Since $\varphi$ and $\ev_e\varphi$ both are $A$-algebra
 homomorphism we obtain, for all $r$ and $r'$ in $R$, that
 $\varphi(rr') =\varphi(r)\varphi(r')$, respectively $\varphi(rr')(e)
 =\varphi(r)(e)\varphi(r')(e)$, such that
 $$\varphi(r)\varphi(r')(e)
 =\varphi(r)(e)\varphi(r')(e).\tag{\sectno.8.1}$$ 
 Moreover we have, since $\ev_e\varphi$ is surjective by assumption,
 that for $x\in F$ and $v\in \End_A(F)$ there are elements $r_v$
 and $r_x$ in $R$ such that $x=\ev_e\varphi(r_x)= \varphi(r_x)(e)$,
 respectively that $v(e) =\ev_e\varphi(r_v) =\varphi(r_v)(e)$.

 Assume first that $v\varphi(r) =\varphi(r)v$ for all $r\in R$. Then
 $v(x) =v(\varphi(r_x)(e)) =\varphi(r_x)v(e)
 =\varphi(r_x)\varphi(r_v)(e)$. It follows from \pref{\sectno.8.1}
 that $v(x) =\varphi(r_x)(e)\varphi(r_v)(e) =xv(e)$, that we wanted to
 show.

 Assume next that $v(x) =v(e)x$ for all $x\in F$. For all $r\in R$ we
 obtain $v\varphi(r)(x) =v(e)\varphi(r)(x)
 =\varphi(r_v)(e)\varphi(r)(\varphi(r_x)(e))
 =\varphi(r_v)(e)\varphi(rr_x)(e)$. It follows from \pref{\sectno.8.1}
 that $v\varphi(r)(x) =\varphi(r_v)\varphi(rr_x)(e)
 =\varphi(r_vrr_x)(e)$. On the other hand we have $\varphi(r)v(x)
 =\varphi(r)v(e)x =\varphi(r)(\varphi(r_e)(e)\varphi(r_x)(e))$. It
 follows from \pref{\sectno.8.1} that $\varphi(r)v(x)
 =\varphi(r)(\varphi(r_e)\varphi(r_x)(e)) =\varphi(rr_vr_x)(e)$. We
 have thus shown that $v\varphi(r)=\varphi(r)v$.

 The last part of the proposition follows since $e$ is part of a basis
 for $F$ and thus all the elements in $F$ are on the form $v(e)$ for
 some $v\in \End_A(F)$.
 \enddemo

 \proclaim{\sectno.9 Corollary}Let $R[Z]$ be the polynomial ring in
 the variable $Z$ over $R$, and let $\beta:F\to R$ be as in
 \cref{\sectno.5}. With the notation of Proposition
 \cref{\sectno.7} we have a natural isomorphism of functors
  $$\Cal H\Cal E_{R[Z]}^\beta\to \Cal H\Cal
  E_R^\beta\times \Cal H_{\Sym_A(F\dual)}\tag{\sectno.9.1}$$
 that for each $A$-algebra $B$ maps $\psi:B\otimes_A R[Z]\to
  \End_B(B\otimes_AF)$ to the pair $(\varphi,\chi)$ where
  $\varphi=\psi\vert B\otimes_AR$ and where $\chi$ is determined as
  follows:

 Let $z\in B\otimes_AF$ be the element determined by
  $\psi(Z)=\chi_\varphi(z)$, where $\chi_\varphi$ is the injective
  $B$-algebra homomorphism of Proposition \cref{\sectno.8} for the
  algebra $B$. Moreover let $u_z:B\to B\otimes_AF$ be the $B$-module
  homomorphism determined by $u_z(1)=z$, and let
  $u_z\dual:B\otimes_AF\dual \to B$ be its dual. Then $\chi$ is
  determined by $\chi(w) =u_z\dual(1_B\otimes w)$ for all $w\in
  F\dual$.

 In particular, if $H_R^\beta$ represents the functor $\Cal H\Cal
 E_R^\beta$ then $H_R^\beta \otimes_A\Sym_A(F\dual)$ represents the
 functor $\Cal H\Cal E_{R[Z]}^\beta$, and the universal families over
 $H_R^\beta$ and $H_R^\beta\otimes_A\Sym_A(F\dual)$ are related via
 the correspondence \pref{\sectno.9.1}.  \endproclaim

 \demo{Proof}It follows from Proposition \cref{\sectno.8} that
 $\psi:B\otimes_AR[Z] \to \End_B(B\otimes_AF)$ belongs to $\Cal H\Cal
 E_{R[Z]}^\beta(B)$ if and only if $\varphi=\psi\vert B\otimes_AR$
 belongs to $\Cal H\Cal E_R^\beta(B)$ and $\psi(Z)$ is contained in
 the image of the $B$-algebra homomorphism
 $\chi_\varphi:B\otimes_AF\to \End_B(B\otimes_AF)$ of Proposition
 \cref{\sectno.8}. The latter condition is the same as saying that
 $\psi(Z)$ is determined by an element in $B\otimes_AF$. Such an
 element is defined uniquely by a $B$-module homomorphism $B\to
 B\otimes_AF$, or its dual $B\otimes_AF\dual \to B$ as in the
 corollary, and the dual corresponds to an $A$-module homomorphism
 $F\dual\to B$. Finally such an $A$-module homomorphism corresponds to
 an $A$-algebra homomorphism $\Sym_A(F\dual)\to B$.  \enddemo

 \def\sectno{2}

 \head \sectno. Representing functors of maps to endomorphisms\endhead 

 In this section we show how the $A$-module homomorphisms $M\to
 \End_A(F)$ from an $A$-module $M$ to the endomorphisms of a free
 $A$-module $F$ can be described by $A$-algebra homomorphisms
 $\Sym_A(M\otimes_A\End_A(F)\dual)\to A$.

 \definition{\sectno.1 Notation}Let $M$ be a fixed $A$-module. We
 write $N\dual =\Hom_A(N,A)$ for the dual of an $A$-module $N$, and
 let $\Sym_A(N)$ be the symmetric algebra of $N$ over $A$. We shall
 consider $N$ as the submodule of the graded $A$-algebra $\Sym_A(N)$
 consisting of elements of degree one.  \enddefinition

 \definition{\sectno.2 Canonical isomorphisms}
 Since $\End_A(F)$ is a free $A$-module of finite rank the
 \dfn{evaluation map} 
 $$\ev: \End_A(F)\otimes_A\End_A(F)\dual \to A$$
 that maps $u\otimes\varphi$ to $\varphi(u)$ corresponds, by duality,
 to a  natural  $A$-module homomorphism
 $$t:A\to \End_A(F)\dual\otimes_A\End_A(F).\tag{\sectno.2.1}$$
 \enddefinition

 \proclaim{\sectno.3 Lemma}For every $A$-algebra $B$ the
 homomorphism of  $B$-modules
 $$\Hom_A(M\otimes_A\End_A(F)\dual, B) \to
 \Hom_A(M,B\otimes_A\End_A(F))\tag{\sectno.3.1}$$
  that maps $u:M\otimes_A\End_A(F)\dual\to B$ to the composite
 homomorphism 
 $$M\risom M\otimes_AA @>{\id_M\otimes t}>>
 M\otimes_A\End_A(F)\dual\otimes_A\End_A(F) 
 @>{u\otimes\id_{\End_A(F)}}>> B\otimes_A\End_A(F)$$
 is an isomorphism. 
 \endproclaim

 \demo{Proof}Choosing a basis for $\End(F)$, and the dual basis for
 $\End(F)\dual$ it is easy to check that the homomorphism of the lemma
 has as inverse the isomorphism that maps an $A$-module homomorphism
 $v:M\to B\otimes_A\End_A(F)$ to
 $(\id_B\otimes\ev)(v\otimes\id_{\End(F)\dual})$. 
 \enddemo

 \definition{\sectno.4 Universal homomorphism}
  Recall that we identify  $M\otimes_A\End(F)\dual$ with a submodule
  of  $\Sym_A(M\otimes_A\End_A(F)\dual)$. Hence we obtain from the
  $A$-module   homomorphism 
 $$M\otimes_AA @>{\id_M\otimes t}>>
 M\otimes_A\End_A(F)\dual\otimes_A\End_A(F).\tag{\sectno.4.1}$$ 
 an $A$-module
 homomorphism 
 $M\to \Sym_A(M\otimes_A\End_A(F)\dual)\otimes_A\End_A(F)$,
 and consequently a natural homomorphism of
 $\Sym_A(M\otimes_A\End_A(F)\dual)$-modules
 $$
 \mu:\Sym_A(M\otimes_A\End_A(F)\dual)\otimes_AM 
 \to  \Sym_A(M\otimes_A \End_A(F)\dual)\otimes_A\End_A(F)
 $$ 
  uniquely determined by
  $\mu(1\otimes x) =x\otimes t(1)$.
   \enddefinition
 
 \proclaim{\sectno.5 Proposition}There is a natural isomorphism of
 functors
 $$\Cal H_{\Sym_A(M\otimes_A\End_A(F)\dual)} \to \Cal
 Hom(M,\End(F))\tag{\sectno.5.1}$$
 such that, for every $A$-algebra $B$, a $B$-algebra homomorphism
 $\varphi:\Sym_A(M\otimes_A\End_A(F)\dual)\to B$  is mapped to the
 $B$-module  homomorphism $\Cal Hom_\varphi(\mu):B\otimes_AM \to
 B\otimes_A\End_A(F)$. 
 \endproclaim

 \demo{Proof}We shall show that the homomorphism \pref{\sectno.5.1} is
 the composite of the following three isomorphisms:
 \roster
 \item $\Hom_{A\-alg}(\Sym_A(M\otimes_A\End_A(F)\dual),B)\to
 \Hom_{A}(M\otimes_A\End_A(F)\dual,B)$ that we obtain from the
 definition of symmetric algebras.
 \item The map \pref{\sectno.3.1}.
  \item The canonical isomorphism $\Hom_A(M,B\otimes_A\End_A(F))\to
 \Hom_B(B\otimes_AM, B\otimes_A\End_A(F))$. 
 \endroster
 Let $u=\varphi\vert (M\otimes_A\End(F)\dual)$. The image of
 $\varphi$ by the composite map of the three isomorphisms is the
 $B$-module homomorphism determined on the element $1_B\otimes x$ by
 $(u\otimes\id_{\End_A(F)})(\id_M\otimes t) (x\otimes 1_A)$ for all
 $x\in M$, and this homomorphism
 is equal to $\Cal Hom_\varphi(\mu)$ since it follows from
 \pref{\sectno.4.1} that for all $x\in M$ we have
 $\Cal Hom_\varphi(\mu)(1_B\otimes x) =(\varphi\otimes
 \id_{\End_A(F)})\mu(1_{\Sym_A(M\otimes_A\End_A(F)\dual)} \otimes x)
 =(\varphi\otimes\id_{\End_A(F)})(x\otimes t(1_A))$.

 As the three isomorphisms are functorial in $B$ we have proved the
  proposition. 
  \enddemo

 \def\sectno{3}

  \head \sectno.  Representing  functors of maps to commuting
  endomorphisms\endhead  

 In section \cref{2} we gave the connection between $A$-module
 homomorphisms $M\to \End_A(F)$ and $A$-algebra homomorphisms
 $\Sym_A(M\otimes_A\End_A(F)\dual)\to A$. Here we show how the
 $A$-module homomorphisms $M\to \End_A(F)$ such that the image
 consists of commuting matrices correspond to $A$-algebra
 homomorphisms $H\to A$, where $H$ is a natural residue algebra of
 $\Sym_A(M\otimes_A\End_A(F)\dual)$ . Hence we obtain a correspondence
 between $A$-algebra homomorphisms $\Sym_A(M)\to \End_A(F)$ and
 $A$-algebra homomorphisms $H\to A$.

  \subhead{\sectno.1 The ideal of zeroes of a homomorphism}\endsubhead
 Let $u:N\to P$ be a homomorphism of an $A$-module $N$ into a free
 $A$-module $P$ of finite rank.  We denote by $\frak
 I_Z(u)$ the ideal in $A$ where $u$ is zero. More precisely, the ideal
 $\frak I_Z(u)$ is the image of the composite homomorphism
 $N\otimes_A P\dual
 @>{\id_{P\dual}\otimes u}>> P \otimes_A P\dual @>{\ev}>> A$. Then,
 for every $A$-algebra $\varphi:A \to  B$, the homomorphism
  $$B\otimes_A N @>{\id_B\otimes u}>> B\otimes_A P$$ is zero if and
 only if $A @>{\varphi}>> B$ factors via the residue map $A\to A/\frak
 I_Z(u)$.

 \definition{\sectno.2 Definition}Let $H$ be the residue algebra of
 $\Sym_A(M\otimes_A\End_A(F)\dual)$ modulo the smallest ideal
 containing the ideals
 $$\frak I_Z(\mu(1\otimes x)\mu(1\otimes y) -\mu(1\otimes
 y)\mu(1\otimes x))\quad\text{for all}\quad x,y\quad\text{in}\quad
 M,$$
 where $\mu$ is the homomorphism of \cref{2.4}. Moreover, let 
 $$\rho_H:\Sym_A(M\otimes_A\End_A(F)\dual)\to H$$
 be the residue class homomorphism. Denote by 
 $$\mu_H : H\otimes_A \Sym_A(M) \to H\otimes_A\End_A(F)$$
 the $H$-algebra homomorphism uniquely defined by
 $$\mu_H(1\otimes x) =\Cal Hom_{\rho_H}(\mu)(1\otimes x)\quad\text{for
 all}\quad x\in M.$$
 \enddefinition
 
 \proclaim{\sectno.3 Proposition}We have a natural isomorphism of
 functors
 $$\Cal H_H\to \Cal Hom_{\alg}(\Sym_A(M), \End(F))$$
 that for every $A$-algebra $B$ maps an $A$-algebra homomorphism
 $\varphi:H\to B$ to the $B$-algebra homomorphism
 $\Cal Hom_\varphi(\mu_H): B\otimes_A\Sym_A(M) \to
 B\otimes_A\End_A(F)$. 
 \endproclaim 

 \demo{Proof}It follows from the definition of $H$ and Proposition
 \cref{2.5} that there is a bijection between the set
 $\Hom_{A{\-alg}}(H,B)$ of $A$-algebra homomorphisms $\varphi:H\to B$
 and the set of $B$-module homomorphisms $u:B\otimes_AM \to
 B\otimes_A\End(F)$ such that the elements $u(1\otimes x)$, for all
 $x$ in $M$, commute. Under this bijection $\varphi$ corresponds to
 the homomorphism $u$ given by $u(1\otimes x) =(\Cal
 Hom_{\varphi\rho_H} (\mu))(1\otimes x) =(\Cal Hom_\varphi\Cal
 Hom_{\rho_H}(\mu))(1\otimes x) =(\Cal Hom_\varphi (\mu_H)) (1\otimes
 x)$.

 From the definition of symmetric products it follows that the set of
 $B$-module homomorphisms $u:B\otimes_AM \to B\otimes_A\End_A(M)$ such
 that the elements $u(1\otimes x)$ commute for all $x$ in $M$
 corresponds bijectively to the set of $B$-algebra homomorphisms
 $\psi: B\otimes_A\Sym_A(M) \to B\otimes_A\End_A(F)$ such that
 $\psi\vert (B\otimes_AM) =u$. We have thus proved that the
 homomorphism $\Cal H_H(B)\to \Cal Hom_{B\-alg}(\Sym_A(M),\End(F))$
 described in the the proposition is an isomorphism. It is clear from
 the construction of the homomorphism that it is functorial in $B$.
 \enddemo

 \def\sectno{4}

  \head \sectno. Sections and closed subschemes\endhead

 In section \cref{3} we described the connection between $A$-algebra
 homomorphisms $\Sym_A(M)\to \End_A(F)$ and $A$-algebra homomorphisms
 $H\to A$. Here we show how we, for any residue algebra $R$ of
 $\Sym_A(M)$, can construct a natural residue algebra $H_R$ of $H$
 such there is a similar connection between $A$-algebra homomorphisms
 $R\to \End_A(F)$ and $A$-algebra homomorphisms $H_R\to A$. This we
 refine further so that we, for every $A$-module homomorphism
 $\beta:F\to R$, obtain a correspondence between $A$-algebra
 homomorphisms $\varphi:R\to \End_A(F)$ such that
 $\ev_e\varphi\beta=\id_F$ and $A$-algebra homomorphisms $H_R^\beta\to
 A$ from a natural residue algebra $H_R^\beta$ of $H_R$.

 \definition{\sectno.1 Definition}Let $\frak I$ be an ideal in
 $\Sym_A(M)$ and let $\iota:\frak I\to \Sym_A(M)$ be the homomorphism
 given by the inclusion of $\frak I$ in $\Sym_A(M)$.  Write
 $R=\Sym_A(M)/\frak I$ and  denote by  $v$  the  composite  of the
 $H$-module  homomorphisms 
 $$ H\otimes_A\frak I @>{\id_{H}\otimes\iota}>>
 H\otimes_A \Sym_A(M)  @>{\mu_{H}}>> 
 H \otimes_A \End_A(F),$$
 where $H$ is defined in \cref{3.2}. 
  Let $H_R$ be the residue ring of $H$ modulo the
 ideal  $\frak I_Z(v)$, of zeroes of $v$, and let 
 $$\rho_{H_R} :H\to H_R$$
 be the residue homomorphism. Then the composite homomorphism
 $$H_R\otimes_A\frak I @>{\id_{H_R}\otimes\iota}>> H_R\otimes_A
 \Sym_A(M) @>{\Cal Hom_{\rho_{H_R}}(\mu_H)}>> H_R\otimes_A\End_A(F)$$
 is zero and thus induces an $H_R$-algebra homomorphism
 $$\mu_{H_R}:H_R\otimes_AR\to H_R\otimes_A\End_A(F).$$ 
 \enddefinition

 \proclaim{\sectno.2 Proposition}We have a natural isomorphism of
 functors
 $$\Cal H_{H_R}\to \Cal Hom_{\alg}(R,\End(F))$$
 such that for every $A$-algebra $B$ the image of an $A$-algebra
 homomorphism $\varphi:H_R\to B$ is
 $\Cal Hom_\varphi(\mu_{H_R}):B\otimes_AR\to B\otimes_A\End_A(F).$
  \endproclaim

 \demo{Proof}By the definition of $H_R$ an $A$-algebra homomorphism
 $\psi:H\to B$ factors via the residue homomorphism
 $\rho_{H_R} :H\to H_R$  if and only if the composite homomorphism
 $$B\otimes\frak I @>{\id_ B\otimes \iota}>> B\otimes_A\Sym_A(M)
 @>{\Cal Hom_\psi(\mu_{H})}>> B\otimes_A\End_A(F)$$ 
 is zero, that is, if and only 
 if the  homomorphism $\Cal Hom_\psi(\mu_{H}):
 B\otimes_A\Sym_A(M)\to B\otimes_A\End_A(F)$ 
 factors via  the residue map $B\otimes_A\Sym_A(M)\to
 B\otimes_A(\Sym_A(M)/\frak I)$.  The proposition
 thus follows from  Proposition  \cref{3.3}. \enddemo

 \definition{\sectno.3 Definition} Recall that $F$ is a free
 $A$-module of finite rank with a distinguished element $e$ that is
 part of a basis. We fix an $A$-module homomorphism
 $$\beta:F\to R\quad\text{such that}\quad \beta(e) =1_R.$$
   Let  $u$ be the composite of the $H_R$-module homomorphisms 
 $$ H_R\otimes_AF@>{\id_{H_R}\otimes\beta}>>
 H_R \otimes_A R@>{\mu_{H_R}}>>  H_R\otimes_A\End_A(F)
 @>{\id_{H_R}\otimes \ev_e}>> H_R\otimes_AF.$$
 We denote by $H_R^\beta$ the residue ring of $H_R$ modulo the ideal
 $\frak I_Z(\id_{H_R\otimes_AF}- u)$ and let 
 $$\rho_{H_R^\beta} :H_R \to H_R^\beta$$
 be the residue  homomorphism. Moreover, we write
 $$\mu_{H_R^\beta} =\Cal Hom_{\rho_{H_R^\beta}}(\mu_{H_R}):
  H_R^\beta\otimes_AR \to H_R^\beta\otimes_A\End_A(F).$$
  \enddefinition

  \proclaim{\sectno.4 Proposition}We have a natural isomorphism of
  functors
 $$\Cal H_{H_R^\beta}\to \Cal Hom_{\alg}^\beta(R ,\End_A(F))$$
 that for every $A$-algebra $B$ maps an $A$-algebra homomorphism
 $\varphi:H_R^\beta\to B$ to $\Cal
 Hom_\varphi(\mu_{H_R^\beta}):B\otimes_AR\to B\otimes_A\End_A(F)$.
 \endproclaim

 \demo{Proof}By the definition of $H_R^\beta$ an $A$-algebra
 homomorphism  $\psi:H_R\to B$ factors via the residue homomorphism
 $\rho_{H_R^\beta} :H_R \to H_R^\beta$ in a homomorphism
 $\varphi:H_R^\beta\to B$ if and only if the 
 composite of the $B$-module homomorphisms
 $$B\otimes_AF@>{\id_{B}\otimes\beta}>>
 B\otimes_A R @>{\Cal Hom_\psi(\mu_{H_R})}>>
 B\otimes_A\End_A(F)  @>{1_{B}\otimes \ev_e}>> B\otimes_AF$$
 is the identity and $\Cal Hom_\psi (\mu_{H_R})= \Cal
 Hom_{\varphi\rho_{H_R^\beta}}(\mu_{H_R}) =\Cal Hom_\varphi
 (\mu_{H_R^\beta})$. The proof thus follows  from Proposition
 \cref{\sectno.2}. 
 \enddemo

 We sum up what we have done in the following result.

 \proclaim{\sectno.5 Theorem}Let $F$ and $M$ be $A$-modules with $F$
 free of finte rank, with a distinguished element $e$ that is part of a
 basis. Moreover, let $\frak I$ be an ideal in $\Sym_A(M)$ and let
 $R=\Sym_A(M)/\frak I$. For all $A$-module homomorphism $\beta:F\to R$
 such that $\beta(e) =1_R$ we have a natural isomorphism of functors
 $$\Cal H_{H_R^\beta} \to \Cal Hom^\beta(R,F)$$
 that is determined by mapping $\id_{H_R^\beta}$ to the homomorphism
 $H_R^\beta\otimes_AR\to H_R^\beta\otimes_AF$ obtained from
 $(\id_{H_R^\beta}\otimes\ev_e) \Cal Hom_\rho(\mu):H_R^\beta\otimes
  M\to H_R^\beta\otimes_AF$, where we let
 $\rho={\rho_{H_R^\beta}\rho_{H^\beta} 
 \rho_H}:\Sym_A(M\otimes_A\End_A(F)\dual)\to H_R^\beta$ denote the
 residue homomorphism. That is, the
 functor $\Cal Hom^\beta(R,F)$ is represented by the $A$-algebra
 $H_R^\beta$ and the {\it universal family} is is the homomorphism
 $H_R^\beta\otimes_AR \to H_R^\beta\otimes_AF$ obtained from
 $(\id_{H_R^\beta}\otimes \ev_e)\Cal Hom_\rho
 (\mu):H_R^\beta\otimes_AM \to H_R^\beta\otimes_AF$.
  \endproclaim

 \demo{Proof}It follows from Proposition \cref{1.7} that we have a
 natural isomorphism of functors $\Cal H\Cal E_R^\beta \to \Cal
 Hom^\beta(R,F)$ and from Proposition \cref{4.4} we have a natural
 isomorphism of functors $\Cal H_{H_R^\beta} \to \Cal Hom_{\alg}^\beta
 (R,\End(F))$. Consequently the theorem follows from the canonical
 isomorphism of functors $\Cal Hom_{\alg}^\beta(R,\End(F)) \to \Cal
 H\Cal E_R^\beta$ that we obtain from the natural isomorphism of
 $B$-algebras $B\otimes_A\End_A(F)\to\End_B(B\otimes_AF)$.  \enddemo

\def\sectno{5}

 \head \sectno. The Hilbert Functor \endhead

 In this section we show how our results for the Hilbert scheme of $n$
 points in $\Spec(R)$ can be used to obtain the Hilbert scheme of any
 projective scheme over an arbitrary base.

 \subhead{\sectno.1 The Hilbert functor in the affine case}
\endsubhead 
 Let $R$ be an $A$-algebra. We let $\Cal Hilb_{R/A}^n$ be the functor
 from $A$-algebras to sets, that to an $A$-algebra $B$ associates the
 set $\Cal Hilb_{R/A}^n(B)$ of surjective $B$-algebra homomorphism
 $\varphi:B\otimes_AR\to Q$, where $Q$ is a locally free $B$-module of
 rank  $n$.  For every $A$-module homomorphism
 $$\beta:F\to R\quad\text{such that}\quad \beta(e) =1_R$$ 
 we let $\Cal Hilb_{R/A}^\beta$ be the subfunctor of  $\Cal
 Hilb_{R/A}^n$ that to $B$  associates the set $\Cal
 Hilb_{R/A}^\beta(B)$ of those $\varphi$ such 
 that
  $$B\otimes_A F @>{\id_{B}\otimes_A \beta}>> B\otimes_A R
 @>{\varphi}>> Q$$
 is surjective, and thus an isomorphism. 

 \subhead{\sectno.2 The existence of sections}\endsubhead
 For every homomorphism $\varphi:B\otimes_A R\to Q$ in
 $\Hilb_{R/A}^n(B)$ and every maximal ideal in $B$ we can find an
 element $f\in B$, not in the maximal ideal, and an $A$-module
 homomorphism $\beta_Q:F\to R$ with $\beta_Q(e) =1$, such that the 
 composite homeomorphism
 $$B_f\otimes_A F @>{\id_{B_f}\otimes_A \beta_Q}>> B_f\otimes_A R
 @>{\varphi_f}>> Q_f$$ 
 of $B_f$-modules is surjective, and thus an isomorphism.  Hence the
 functors $\Hilb_{R/A}^\beta$, for all choises of $\beta$, form an
 open cover of 
 $\Hilb_{R/A}^n$. We can even choose the $\beta$ 
 to map a fixed basis of $F$ to a given set of generators of $R$ to
 obtain an open covering, and when $R=\polring$ we can map the basis
 of $F$  to monomials in $X_1, \dots, X_m$ of degree strictly less
 than $n$. 

 \proclaim{\sectno.3 Lemma}We have a natural isomorphism of functors 
 $$\Cal Hom^\beta(R,F) \to \Cal Hilb_{R/A}^\beta$$
 given, for every $A$-algebra $B$,
 by the canonical homomorphism $\Cal Hom_B^\beta(R,F) \to \Cal
 Hilb_R^\beta(B)$ such that the image of a $B$-module homomorphism
 $u:B\otimes_A R\to B\otimes_A F$ is equal to itself when
 $B\otimes_AF$ is given the unique $B$-algebra structure such that $u$
 becomes a $B$-algebra homomorphism.
 \endproclaim

  \demo{Proof} To see that the morphism given in the lemma is an
 isomorphism we construct an inverse. For every $B$-algebra
 homomorphism $\varphi:B\otimes_A R\to Q$ in $\Hilb_R^\beta(B)$ the
 composite homomorphism $B\otimes_A F@>{\id_B\otimes \beta}>>
 B\otimes_A R @>{\varphi}>> Q$ is an isomorphism. Via this isomorphism
 $B\otimes_A F$ obtains a unique $B$-algebra structure such that
 $B\otimes_A R @>{(\varphi(\id_B\otimes \beta))^{-1}\varphi}>>
 B\otimes_A F$ is in $\Cal Hom_B^\beta(R,F)$. It is clear that this
 defines an inverse to the map of the lemma.\enddemo

 \proclaim{\sectno.4 Theorem}Let $R$ be an $A$-algebra. The functor
 $\Cal Hilb_{R/A}^n$ is representable. More precisely, the functor
 $\Cal Hilb_{R/A}^n$ is covered by the open subfunctors $\Cal
 Hilb_{R/A}^\beta$ for all $A$-module homomorphism $\beta:F\to R$, and
 for each $\beta$ we have a natural isomorphism of functors
 $$\Cal H_{H_R^\beta} \to \Cal Hilb_{R/A}^\beta$$
 that maps $\id_{H_R^\beta}$ to $\mu_{H_R^\beta}: H_R^\beta\otimes_AR
 \to H_R^\beta\otimes_AF$. That is, the functor $\Cal
 Hilb_{R/A}^\beta$ is represented by the $A$-algebra $H_R^\beta$ and
 the universal homomorphism is $\mu_{H_R^\beta}$.
 \endproclaim

 \demo{Proof}The first part of the theorem we observed in Section
 \cref{\sectno.2}, and the isomorphism of the theorem is the composite
 of the natural isomorphism of functors of Theorem \cref{4.5} and
 Lemma \cref{\sectno.3}.  \enddemo

 \subhead{\sectno.5 The Hilbert functor}\endsubhead Let $X$ be a
 scheme over a base scheme $S$ and let $p:X\to S$ be the structure
 homomorphism. For a morphism $T\to S$ from a scheme $T$ we let
 $p_T:T\times_SX\to T$ denote the projection. The $T$-points $\Cal
 Hilb_{X/S}^n(T)$ of the \dfn{Hilbert functor} $\Cal Hilb_{X/S}^n$ of
 $n$ points of $X$ over $S$ consists of the closed subschemes $Z$ of
 $T\times_SX$ such that $Z$ is finite over $T$ and $(p_T)_\ast(\Cal
 O_Z)$ is a locally free $\Cal O_T$-module of rank $n$ (see
 \cite{LST}). 

 \proclaim{\sectno.6 Lemma}Let $R$ be a graded $A$-algebra. For every
 prime ideal $\frak p$ of $A$  write $\boldkappa(\frak p) =A_{\frak
 p}/\frak pA_{\frak p}$. 

 Let $Z$ be a closed subscheme of $\Proj(\boldkappa(\frak
 p)\otimes_AR)$ that is finite over $\Spec(\boldkappa(\frak p))$. Then
 there is an element $a\in A$ not in $\frak p$ and an element $f\in
 R_a$ such that $Z$ is contained in the open subscheme
 $\Spec(\boldkappa(\frak p)\otimes_{A_a} (R_a)_{(f)})
 =\Spec(\boldkappa(\frak p))\times_{\Spec(A_a)}\Spec((R_a)_{(f)})$ of
 $\Proj(\boldkappa(\frak p)\otimes_{A_a} R_a) =\Spec(\boldkappa(\frak
 p))\times_{\Spec(A_a)} \Proj(R_a)$.  \endproclaim

 \demo{Proof}Since $Z$ is finite over $\Spec(\boldkappa(\frak p))$ the
 fiber of the induced morhpism $Z\to \Spec(\boldkappa(\frak p))$
 consists of a finite number of points, corresponding to homogeneous
 prime ideals $\frak q_1, \dots, \frak q_k$ in $\boldkappa(\frak
 p)\otimes_AR$ that do not contain the irrelevant ideal. Their union
 consequently do not contain the irrelevant ideal. Hence we can find a
 homogeneous element $g\in \boldkappa(\frak p)\otimes_AR$ of positive
 degree that is not contained in any of the ideals $\frak q_1,\dots,
 \frak q_k$. Thus $Z$ is contained in the open subscheme
 $\Spec(\boldkappa(\frak p)\otimes_AR)_{(g)})$ of
 $\Proj(\boldkappa(\frak p)\otimes_AR)$. 

 Clearly we can find an element $a\in A$ not in $\frak p$ and an
 element $f\in R_a$ such that $1_{\text{\miibold \char'024}(\frak
 p)}\otimes f$ is the image of $g$ by the natural isomorphism
 $\boldkappa(\frak p)\otimes_AR\to \boldkappa(\frak
 p)\otimes_{A_a}R_a$. However, then $\Spec((\boldkappa(\frak
 p)\otimes_AR)_{(g)}) =\Spec(\boldkappa(\frak
 p)\otimes_{A_a}(R_a)_{(f)})$, and we have proved the lemma.  \enddemo

 \proclaim{\sectno.7 Theorem}Let $S$ be a scheme and $\Cal R$ a
 quasi-coherent graded $\Cal O_S$-algebra. Then the functor $\Cal
 Hilb_{\Proj(\Cal R)/S}^n$ is representable.

 More precisely, the functor $\Cal Hilb_{\Proj(\Cal R)/S}^n$ is
 covered by the representable open subfunctors $\Cal
 Hilb_{\Spec(R_{(f)})/\Spec(A)}^n$, where $\Spec(A)$ is an open affine
 subscheme of $S$, where $R=\Gamma(\Spec(A),\Cal R\vert \Spec(A))$,
 and where $f$ is a homogeneous element of $R$.
 \endproclaim

 \demo{Proof} For every affine open subset $U$ of $S$ we have that
 $\Cal Hilb_{p^{-1}(U)/U}^n$ is an open subfunctor of $\Hilb_{X/S}^n$
 and these subfunctors, for all $U$ in an open covering of $S$, cover
 $\Hilb_{X/S}^n$. In order to represent $\Cal Hilb_{X/S}^n$ we can
 thus assume that $S=\Spec(A)$ is affine. Write
 $R=\Gamma(\Spec(A),\Cal R)$ such that $X=\Proj(R)$. For every $a\in
 A$ and every $f\in R_a$ the functor $\Cal
 Hilb_{\Spec((R_a)_{(f)})/\Spec(A_a)}^n$ is an open subfunctor of
 $\Cal Hilb_{X/S}^n$. It follows from Theorem \cref{\sectno.4} that in
 order to prove the theorem it suffices to prove that these
 subfunctors cover the functor $\Cal Hilb_{X/S}^n$. In order to show
 this it suffices to show that for every prime ideal $\frak p$ in $A$
 and every closed subscheme $Z$ of $\Proj(\boldkappa(\frak
 p)\otimes_AR)$ that is finite over $\Spec(\boldkappa(\frak p))$, and
 with $\Gamma(Z,\Cal O_Z)$ of dimension $n$ over $\boldkappa(\frak
 p)$, there is an $a\in A$, not in $\frak p$ and an $f\in R_a$ such
 that $Z$ is contained in the open subscheme $\Spec(\boldkappa(\frak
 p))\times_{\Spec(A_a)}\Spec((R_a)_{(f)}) =\Spec(\boldkappa(\frak
 p)\otimes_{A_a}(R_a)_{(f)})$ of $\Proj(\boldkappa(\frak
 p)\otimes_{A_a} R_a)$. It consequently follows from Lemma
 \cref{\sectno.6} that the functors $\Cal
 Hilb_{\Spec(R_a)_{(f)}/\Spec(A_a)}^n$ cover $\Cal Hilb_{X/S}^n$, and
 we have proved the theorem.
 \enddemo

\def\sectno{6}

 \head \sectno. The Hilbert scheme in coordinates\endhead

 In this section we express our construction of the Hilbert scheme of
 $n$ points in $\Spec(R)$ in local coordinates. We obtain explicit
 expression of an open affine covering of the Hilbert scheme in terms
 of variables and relations.

 \definition{\sectno.1 Notation}We choose a basis $e=T_1, T_2,\dots,
 T_n$ of the $A$-module $F$, and let $T_1\dual, \dots, T_n\dual$ be
 the dual basis of $F\dual$. Correspondingly we obtain a basis
 $T_{ij}$ for $i,j=1,\dots,n$ of $\End_A(F)$, and a dual basis
 $T_{ij}\dual$ for $\End_A(F)\dual$. In particular $\id_{\End_A(F)}=
 \sum_{i=1}^n T_{ii}$ and $\id_{\End_A(F)\dual}=\sum_{i=1}^n
 T_{ii}\dual$. Then
 $$T_{ij}(T_k)=\delta_{jk}T_i\quad\text{and}\quad
 T_{ij}\dual(T_k\dual) =\delta_{ik}T_j\dual.$$
 In these bases we have
 $$t(1_A) =\sum_{i,j=1}^n T_{ij}\dual\otimes T_{ij},$$
 where $t$ is defined in Section \pref{2.2.1}.
 For every $A$-algebra $B$ we consider $T_{ij}$ as a basis for the
 $B$-module $\End_B(B\otimes_A F)$ and thus identify
 $B\otimes_A\End_A(F)$ with 
 $\End_B(B\otimes_AF)$ via the homomorphism that maps $1\otimes
 T_{ij}$ to $T_{ij}$. Then we have for all $b_{ij}$ in $B$ that
  $$\ev_{1_B\otimes T_1}(\sum_{i,j=1}^n b_{ij}T_{ij}) = (\id_B\otimes
 \ev_{T_1}) (\sum_{i,j=1}^n b_{ij}\otimes T_{ij}) =\sum_{i=1}^n
 b_{i1}.$$ 
 Assume that $M$ is a free $A$-module with basis $\{Y_s\}_{s\in S}$
 for some index set $S$. Then $\Sym_A(M\otimes_A\End_A(F)\dual)$ is the
 polynomial ring over $A$ in the independent variables $Y_s\otimes
 T_{ij}\dual$,  for $i,j=1,\dots, n$ and $s\in S$. We write $U_{ij}^s
 =Y_s\otimes T_{ij}\dual$, and $A[U]
 =\Sym_A(M\otimes_A\End_A(F)\dual)$. As above we identify
 $M\otimes_A\End_A(F)\dual$ with the degree $1$ part of
 $\Sym_A(M\otimes_A\End_A(F)\dual)$, and we shall identify $M$ with a
 submodule of $M\otimes_A\End_A(F)\dual$ via the map that takes $x$ to
 $x\otimes\id_{\End_A(F)\dual}$. In  particular  $Y_s = Y_s\otimes
 \id_{\End_A(F)\dual} =Y_s\otimes\sum_{i=1}^n T_{ii}\dual
 =\sum_{i=1}^n U_{ii}^s$, and  $\Sym_A(M)= A[Y]$ is the
 $A$-algebra in  
 $A[U]$ generated by the elements $Y_s$ for $s\in S$. We write
 $$(U_{ij}^s) =\sum_{i,j=1}^n U_{ij}^s \otimes T_{ij} = \left(
 \smallmatrix  U_{11}^s &\cdots &U_{1n}^s \\
 \vdots &\ddots &\vdots\\
 U_{n1}^s &\cdots &U_{nn}^s\\
 \endsmallmatrix \right).$$
 For every polynomial $f(Y)$ in $A[Y]$ we write $f((U_{ij}^s))$ for
 the element in $A[U]\otimes_A\End_A(F)$ obtained by substituting the
 matrix $(U_{ij}^s)$ for the variable $Y_s$. The $A[U]$-module
 homomorphism $\mu:A[U]\otimes_AM \to A[U]\otimes_A\End_A(F)$
 is determined by
 $$\mu(1_{A[U]}\otimes Y_s) =\sum_{i,j=1}^n Y_s\otimes T_{ij}\dual
 \otimes 
 T_{ij}  = \sum_{i,j=1}^n U_{ij}^s\otimes T_{ij} =(U_{ij}^s).$$
 \enddefinition

 \subhead{\sectno.2 Coordinates of the representing
 ring}\endsubhead 
 \roster
 \item Let $\frak I_1$ be the ideal in $A[U]$ generated by the
 coordinates of the matrices
 $$(U_{ij}^s)(U_{ij}^t) - (U_{ij}^t)(U_{ij}^s)$$
 for all $s,t$ in $S$. Then $H=A[U]/\frak I_1$, and
 $\mu_H:H\otimes_AA[Y]\to H\otimes_A\End_A(F)$ is determined by
 $$\mu_H(1_H\otimes f(Y)) =f((U_{ij}^s))$$
 for all $f(Y)$ in $A[Y]$.
 \item Let $\frak I$ be an ideal of $A[Y]$ and let $\frak
 I_2$ be the ideal in $A[U]$  generated by the coordinates of the
 matrices 
 $$f((U_{ij}^s))\quad\text{for}\quad f\in \frak I.$$
 Then $H_{A[Y]} =A[U]/(\frak I_1, \frak I_2)$.
 \item Let $\beta:F\to A[Y]$ be an $A$-module homomorphism and write
 $f_k(Y)= \beta(T_k)$ for $k=1,\dots, n$. Denote by $\frak
 I_3$ the ideal in $A[U]$ generated by the  coefficients of
 $T_1,\dots, T_n$ in 
 $$\ev_{1_{A[U]} \otimes T_1} (f_k((U_{ij}^s)))
 -T_k\quad\text{for}\quad  k=1,\dots, n.$$
 Then $H_{A[Y]}^\beta =A[U]/(\frak I_1, \frak I_2, \frak I_3)$.
 \endroster

 \def\sectno{7}

 \head{\sectno. The generic open subset of the Hilbert scheme of
 affine  space}\endhead

 We use the explicit coordinate description in section \cref{6} of an
 open covering of the Hilbert scheme to give a simple expression of an
 open subset of the \dfn{generic} component of the Hilbert scheme of
 $n$ points in $\Spec(A[Y])$, where $Y_s$ for $s\in S$ is a set of
 independent variables over $A$.

 \definition{\sectno.1 Notation}We shall keep the notation of Section
 \cref{6}, and assume in addition that $S$ has a distinguished element
 $s_1$. Denote by
 $$\beta_1:F\to A[Y]$$
 the $A$-module homomorphism defined by
 $$\beta_1(T_i) =Y_1^{i-1}\quad\text{for}\quad i=1,\dots,n.$$
 We write
 $$C_1=C_1^{s_1} =\left(\smallmatrix
 0&0&0&\cdots&0&0&a_{1}^{s_1}\\ 
 1&0&0&\cdots&0&0&a_{2}^{s_1}\\
 0&1&0&\cdots&0&0&a_{3}^{s_1}\\
 \vdots&\vdots&\vdots&\ddots&\vdots&\vdots&\vdots\\
 0&0&0&\cdots&0&1&a_{n}^{s_1}\\
 \endsmallmatrix\right),$$
 in the basis $T_{ij}$ of $\End_A(F)$. Then
 $$C_1^n =a_1^{s_1}I_n +a_2^{s_1}C_1+\cdots
 +a_n^{s_1}C_1^{n-1}.$$
 Let 
 $$\varphi_1: A[Y_1]\to \End_A(F)$$
 be the $A$-algebra homomorphism determined by
 $\varphi_1(Y_1)=C_1$. Then the composite homomorphism
 $$F@>{\beta_1}>> A[Y_1] @>{\varphi_1}>> \End_A(F) @>{\ev_e}>>F$$
 is the identity on $F$.
 \enddefinition

 \proclaim{\sectno.2 Lemma}Let $\varphi:A[Y]\to \End_A(F)$ be an
 $A$-algebra homomorphism. Then  the composite  homomorphism
 $$F @>{\beta_1}>> A[Y] @>{\varphi}>> \End_A(F) @>{\ev_e}>>
 F\tag{\sectno.2.1}$$
 is the identity on $F$ if and only if $\varphi$ is defined by
 $$\varphi(Y_1)=C_1\quad\text{and}\quad \varphi(Y_s)=a_1^sI_n
 +a_2^sC_1+\cdots +a_n^sC_1^{n-1},\tag{\sectno.2.2}$$
 for all $s\in S\setminus\{s_1\}$, with $a_1^s, \dots, a_n^s$ in $A$.
 \endproclaim

 \demo{Proof}Assume that $\varphi$ is of the form
 \pref{\sectno.2.2}. Since $\varphi(Y_1) =C_1$ the composite of the
 homomorphisms of \pref{\sectno.2.1} is the same as the composite of
 the homomorphisms of \pref{\sectno.1.1}, and thus, as we observed in
 Section \cref{\sectno.1}, equal to the identity of $F$.

 Conversely, assume that the composite of the maps \pref{\sectno.2.1}
 is the identity on $F$. Since $\varphi(Y_s)$ commute with
 $\varphi(Y_1)=C_1$ for all $s\in S$, it follows from Proposition
 \cref{1.8} that $\varphi(Y_s)$ is determined by $\varphi(Y_s)(T_j)
 =\varphi(Y_s)(T_1)T_j$ for $j=1,\dots, n$, where $F$ has the
 $A$-algebra structure determined by the surjection $\ev_e\varphi_1
 =\ev_e(\varphi\vert A[Y_1])$. Since $\varphi$ and $\ev_e\varphi$ are
 $A$-algebra homomorphisms we obtain, as in \pref{1.8.1}, that $T_iT_j
 =\varphi(Y_1^{i-1})(T_1)\varphi(Y_1^{j-1})(T_1)
 =\varphi(Y_1^{i-1})\varphi(Y_1^{j-1})(T_1) =C_1^{i-1}C_1^{j-1}(T_1)
 =C_1^{i-1}(T_j)$ for $i,j=1,\dots, n$ . Wrtie $\varphi(Y_s)(T_1)
 =a_1^sT_1+\cdots+a_1^sT_n$. Then $\varphi(Y_s)(T_1)T_j
 =(a_1^sT_1+\cdots +a_n^sT_n)T_j =(a_a^sI_n+a_2^sC_1+\cdots
 +a_n^sC_1^{n-1})(T_j)$ for $j=1,\dots, n$. Hence $\varphi(Y_s)
 =a_1^sI_n+a_2^sC_1+\cdots+a_n^sC_1^{n-1}$ as we wanted to prove.
 \enddemo

 \proclaim{\sectno.3 Proposition}We have that $H_{A[Y]}^{\beta_1}$ is
 the polynomial ring over $A$ in the independent variables
 $U_{1n}^{s_1}, \dots, U_{nn}^{s_1}$ and $U_{11}^s, \dots, U_{n1}^s$
 for all $s$ in $S$ different from $s_1$. The universal family is
 given by the ideal in $A[Y]\otimes_AH_{A[Y]}^{\beta_1}$ generated by
 the elements $Y_1^n -U_{1n}^{s_1} -U_{2n}^{s_1}Y_1-\cdots
 -U_{nn}^{s_1}Y_1^{n-1}$ and $Y_s -U_{1n}^s -U_{2n}^sY_1-\cdots
 -U_{nn}^sY_1^{n-1}$ for $s\in S\setminus\{s_1\}$.
\endproclaim

 \demo{Proof}The ideal $\frak I_3$ of \cref{6.2}(3) is generated by
 the coefficients of $T_k$ in the polynomials
 $\ev_{1_{A[U]}\otimes_AT_1} ((U_{ij}^s)^{k-1}) -T_k$ for $k=1,\dots,
 n$. It follows from Lemma \cref{\sectno.2} that, modulo the ideal
 $\frak I_3$, we have that $(U_{ij}^{s_1})$ is congruent to the
 \dfn{companion} matrix
   $$C_U =\left(\smallmatrix 0&0&0&\cdots&0&0&U_{1n}^{s_1}\\
 1&0&0&\cdots&0&0&U_{2n}^{s_1}\\
 0&1&0&\cdots&0&0&U_{3n}^{s_1}\\
 \vdots&\vdots&\vdots&\ddots&\vdots&\vdots&\vdots\\
 0&0&0&\cdots&0&1&U_{nn}^{s_1}\\
 \endsmallmatrix\right)$$
 and that 
 $$U_{ij}^s \equiv U_{11}^s I_n +U_{21}^s C_U+\cdots  +U_{n1}^s
 C_U^{n-1} \pmod{\frak I_3}.$$
 In particular the matrices $(U_{ij}^s)$ commute modulo $\frak I_3$,
 that is $\frak I_1\subseteq \frak I_3$, and since $\frak I_2=0$ in
 this case, it follows that $H_{A[Y]}^{\beta_1} =A[U]/\frak I_3$.

 To prove the last part of the Proposition we note that since $I_n,
 C_U, \dots, C_U^{n-1}$ clearly are linearly independent over
 $H_{A[Y]}^{\beta_1}$ the $H_{A[Y]}^{\beta_1}$-algebra homomorphism
 $\varphi_U:A[Y]\otimes_AH_{A[Y]}^{\beta_1}\to
 \End_{H_{A[Y]}}(H_{A[Y]}^{\beta_1}\otimes_A F)$ defined by
 $\varphi_U(Y_1) =C_U$ and $\varphi_U(Y_s) =C_U^s =U_{1n}^sI_n
 +U_{2n}^sC_U+\cdots +U_{nn}^sC_U^{n-1}$ for $s\in S\setminus\{s_1\}$
 has kernel generated by the elements in the last part of the
 proposition.  \enddemo

\def\sectno{8}

  \head{\sectno. The Hilbert scheme of points on a line}\endhead

 To illustrate how easily certain question of Hilbert schemes can be
 handled by our explicit expressions of open coverings of Hilbsert
 schemes we describe the Hilbert scheme of $n$ points in
 $\Spec(S^{-1}A[X])$, where $S$ is a multiplicative set in the
 polynomial ring $A[X]$ in the variable $X$ over $A$.

 \definition{\sectno.1 The open covering} Let $A[X]$ be the polynomial
 ring in the variable $X$ over $A$ and let $S$ be a multiplicatively
 closed subset of $A[X]$ containing $1$. Moreover, let $\{Y_s\}_{s\in
 S}$ be a collection of independent variables with $Y_1=X$. We denote
 by $M$ the free $A$-module generated by the elements $Y_s$ for $s\in
 S$ and, as in \cref{6.2} we write $A[Y] =\Sym_A(M)$. Then there
 is a surjective $A$-algebra homomorphism
 $$\varphi:A[Y] \to S^{-1}A[X]$$
 defined by $\varphi(Y_1) =X$ and $\varphi(Y_s) = 1/s(X)$ for all other
 $s$ in $S$. The kernel of $\varphi$ is generated by the elements
 $$s(X)Y_s -1\quad\text{for all } s\text{ in } S \text{ different from
 }  1.$$
 Let $B$ be an $A$-algebra and let
 $$\psi:B\otimes_AS^{-1}A[X]\to C$$
 be a surjective homomorphism of $B$-algebras, where $C$ is a free
 $B$-module of rank $n$ and write $\psi(1\otimes X)=f$. It follows
 from the  Cayley-Hamilton Theorem that there is a relation $f^n
 +b_1f^{n-1}+\cdots +b_n=0$ in $C$, with $b_1, \dots, b_n$ in
 $B$. Consequently the image of $B\otimes_AA[X]$ by $\psi$ is in the
 $B$-module generated by $1,f,\dots, f^{n-1}$. Let $s$ be in $S$ and
 write $\psi(1\otimes s)=g$. Then $g$ is invertible in $C$. Write
 $\psi(1\otimes s^{-1}) =g^{-1}$. Since $\psi(1\otimes s^{-1})$ also
 satisfies a relations $(g^{-1})^n+a_1 (g^{-1})^{n-1}+\cdots +a_n=0$
 in  $C$,  with $a_1, \dots, a_n$ in $B$, we obtain that
 $\psi(1\otimes s^{-1}) =g^{-1} = -a_1-\cdots -a_n g^{n-1}$, and thus
 the image of $1/s(X)$ lies in the  $B$-module generated by $1,f,\dots,
 f^{n-1}$. When  $F$ is the free $A$-module with basis
 $e=T_1,\dots, T_n$ and
 $$\beta:F\to A[X]$$
 is defined by $\beta(T_i) =X^{i-1}$ for $i=1,\dots, n$, then $\Cal
 Hilb_{S^{-1}A[X]/A}^n =\Cal Hilb_{S^{-1}A[X]/A}^\beta$.
 \enddefinition

 \definition{\sectno.2 Description of the coordinate ring of the
 Hilbert  scheme}
 In this example the ideal $\frak I_3$ in \cref{6.2}(3) is
 generated  by the coefficients of $T_1, \dots, T_n$ 
 in
 $$\ev_{1_{A[U]}\otimes T_1}((U_{ij}^s)^{k-1})
 -T_k\quad\text{for}\quad k=1,\dots, n.$$
 The relation $\ev_{1_{A[U]}\otimes T_1}((U_{ij}^1)^{k-1}) =T_k$
 expressed that the first column in $(U_{ij}^1)^{k-1}$ is equal to the
 column vector with $1$ in the $k$'th coordinate and $0$
 elsewhere. It follows from Lemma \cref{7.2} the matrix $(U_{ij}^1)$ is
 congruent,  modulo the  ideal $\frak I_3$, to the 
 \dfn{companion} matrix 
  $$C_U =\left(\smallmatrix 0&0&0&\cdots&0&0&U_{1n}^1\\
 1&0&0&\cdots&0&0&U_{2n}^1\\
 0&1&0&\cdots&0&0&U_{3n}^1\\
 \vdots&\vdots&\vdots&\ddots&\vdots&\vdots&\vdots\\
 0&0&0&\cdots&0&1&U_{nn}^1\\
 \endsmallmatrix\right).$$
 In this example the elements given in \pref{6.2}(2) are
 $$s((U_{ij}^1))(U_{ij}^s) -1\quad\text{for}\quad s\in S.$$
  Then we have $\det s((U_{ij}^1))\det(U_{ij}^s) \equiv 1
 \pmod{\frak I_2}$. In particular, the elements $d_s
 =\det(s((U_{ij}^1))$ are invertible modulo $\frak I_2$, and
 $(U_{ij}^s) \equiv d_s^{-1}V_s \pmod {\frak I_2}$, where $V_s$ is the
  \dfn{adjoint} matrix of 
 $s((U_{ij}^1))$, and consequently has coordinates that are polynomials
 in $U_{ij}^1$ for $i,j=1,\dots, n$. 
  We see that $A[U]/\frak I_2
 =D^{-1}A[U_{11}^1, \dots, U_{ij}^1, \dots, U_{nn}^1]$ where $D$ is
 the multiplicatively closed subset consisting of all products of the
 $d_s$ for $s\in S$. Moreover we have seen that $A[U]/(\frak I_2,
 \frak I_3) =E^{-1}A[U_{1n}^1, \dots, U_{nn}^1]$, where $E$ consists of
 the the  elements obtined from the elements of $D$ by the
 specialization that takes  $(U_{ij}^1)$ to the companion matrix
 $C_U$. Moreover, since the matrices
 $s((U_{ij}^1))$ commute for $s\in S$, the matrices $(U_{ij}^s)$ 
 commute modulo $\frak I_2$. Consequently $\frak I_1$ of
  \cref{6.2}(2) is contained in $\frak I_2$ and $H_R^\beta
  =A[U]/(\frak I_2, \frak I_3)$.

 In order to get a more attractive presentation of $E^{-1}A[U_{1n}^1,
 \dots, U_{nn}^1]$ we factor the characteristic polyomial of $C_U$ as
 $$p_{C_U}(T) =T^n -U_{nn}^1T^{n-1}- \cdots -U_{1n}^1 =\prod_{i=1}^n
 (T-Z_i)$$
 where $Z_1, \dots, Z_n$ are independent variables over $A$. That is,
 we have  $U_{n-i+1\, n}^1 =(-1)^{i+1}c_i(Z)$ for $i=1,\dots, n$,
 where  $c_i(Z)$ is 
 the $i$'the elementary symmetric function in $Z_1, \dots, Z_n$. It
 follows from the {\it Spectral Mapping Theorem}
 (\cite{EL},\cite{LSvT}, \cite{L} Chapter XIV, \S 3, Theorem 3.10,
 p.~566) that 
 $$\det(s(C_U)) =\prod_{i=1}^n s(Z_i).$$
 It follows that the functor  $\Cal Hilb_{S^{-1}A[X]/A}^n$ is
 represented  by the  polynomial ring 
 $A[c_1(Z), \dots, c_n(Z)]$ localized in the elements
 $\prod_{i=1}^ns(Z_i)$  for all $s\in S$.
 \enddefinition

 \def\sectno{9}

 \head \sectno. A degenerate open subset of the Hilbert scheme of
 affine  space \endhead

 In this section we describe an open subset of a component of the
 Hilbert scheme containing many subschemes of $n$ points with support
 at a fixed point. One of the important features of this set is that
 it can be used (see \cite{I}) to show that the Hilbert scheme of $n$
 points in $\Spec(A[Y_1,\dots, Y_m])$ is reducible when $m\geq 3$ and
 $n$ is large.

 \definition{\sectno.1 Notation}We shall keep the notation of Section
 \cref{6}, and assume in addition that $S=\{1,\dots, m\}$ such that
 $R=A[Y_1,\dots, Y_m]$.  Determine
 the integer $d$ by the inequalities
 $$\binom{d+m-1}{m} <n \leq \binom{d+m}{m}\tag{\sectno.1.1}$$
 and let $s=\binom{d+m}{m} -n$. We order the monomials in the
 variables $Y_1, \dots, Y_m$ lexicographically and let $m_1,
 m_2,\dots$ be the monomials in this ordering. Hence the equalities
 \pref{\sectno.1.1} are equivalent with the condition that the $n$'th
 monomial is of degree $d$.
 \enddefinition

 \subhead{\sectno.2 The section defining the open subset}\endsubhead
 Let
 $$\beta:F\to A[Y]$$
 be the $A$-module homomorphism defined by $\beta(T_i) =m_i$ for
 $i=1,\dots, n$. For each $s\times(\binom{d-1+m}{m-1} -s)$-matrix
 $a=(a_{ij})$ with entries from $A$ we define an $A$-module
 homomorphism 
 $$u_a:A[Y]\to F$$
 by
 $$u_a(m_i) =\cases T_i\quad\text{for}\quad i=1,\dots, n\\
 \sum_{j=\binom{d+m-1}{m} +1}^n a_{ij} T_j\quad\text{for}\quad
 i=n+1,\dots, \binom{d+m}{m}\\
0\quad\text{for}\quad i=\binom{d+m}{m}+1,\dots
 \endcases.$$
 It is clear that the kernel of $u_a$ is the ideal
 generated by the homogeneous monomials of degree strictly greater
 than $d$ and by  the polynomials
 $$m_i-\sum_{j=\binom{d-1+m}{m}+1}^n a_{ij}m_j\quad\text{for}\quad
 i=n+1,\dots, \binom{d+m}{m}.$$ 
 Moreover, it is clear that different matrices $a=(a_{ij})$ give
 different  ideals. Consequently we obtain for each
 $s\times(\binom{d+m-1}{m-1}-s)$-matrix $a=(a_{ij})$ a unique
 $A$-algebra homomorphism
 $$\varphi_a:A[Y]\to \End_A(F)$$ such that
 $\ev_{T_1}\varphi_a=u_a$. In particular we have that $\varphi_a$
 corresponds to an element in $\Cal Hilb_{A[Y]}^\beta(A)$. We see that
 the dimension of $\Cal Hilb_{A[Y]}^\beta$ is at least equal to
 $s(\binom{d+m-1}{m-1}-s)$. In order to get the
 dimension of $\Cal Hilb_{A[Y]}^\beta$ as big as possible we must
 choose $s=\lfloor \binom{d+m-1}{m-1}/2\rfloor$. Easy computations
 (see \cite{I}, \S~3) show that the dimension of $\Cal
 Hilb_{A[Y]}^\beta$ is at least equal to
 $n^{2-(2/m)}(m!/2)^{-(2/m)}(m^2/16)$ when $d\geq 2m^2$.

 \enddocument
 \end

\enddocument